# Novel hybrid function operational matrices of fractional integration: An application for solving multi-order fractional differential equations

Seshu Kumar Damarla[1], Madhusree Kundu[2]

[1]Department of Chemical Engineering, C.V. Raman College of Engineering, Bhubaneswar, Odisha, India.

[1]E-mail: seshukumar.damarla@gmail.com

[2]Department of Chemical Engineering, National Institute of Technology Rourkela, India-769008.

[2]E-mail: mkundu@nitrkl.ac.in

**Abstract**

In the present work, an attempted was made to develop a numerical algorithm by the use of new orthogonal hybrid functions formed from hybrid of piecewise constant orthogonal sample-and-hold functions and piecewise linear orthogonal triangular functions to obtain the numerical solution of multi-order fractional differential equations. The construction of numerical algorithm involves a formula, consisting of generalized one-shot operational matrices, expressing explicitly the Riemann-Liouville fractional order integral of the new orthogonal hybrid functions in terms of new orthogonal hybrid functions themselves. A set of test problems comprising of linear and nonlinear multi-order fractional differential equations with constant and variable coefficients was considered to demonstrate the accuracy and computational efficiency of our algorithm and to compare our results with those acquired by some other well-known methods used for solving multi-order fractional differential equations in the literature. The comparison highlighted that our algorithm exhibits superior performance to those methods.

## 1. Introduction

Since the fractional calculus was born during the conversation held on 30 September 1695 between Leibniz and L'Hospital (Oldham and Spanier 1974), a number of methods has been devised until the present for solving ordinary differential equations involving derivatives and/or integrals of non-integer order. At first analytical methods such as Fourier transforms (Gaul et al. 1989, 1991), modal synthesis and eigenvector expansion (Suarez and Shokooh 1997), Laplace transform method, fractional Green's function, Mellin transform method, power series method (Podlubny 1999), etc. were developed to find the closed form solutions of linear fractional order differential equations. After realizing the fact that the aforesaid analytical methods cannot be applied to get analytical solutions of nonlinear fractional differential equations and even complicated linear fractional order differential equations, numerical techniques like Podlubny's numerical method (Podlubny 1999), Predictor Evaluate–Corrector Evaluate type method (Diethelm et al. 2002), generalized Euler's method (Odibat and Momani 2008) and semi numerical techniques, for instance, Adomian decomposition method (Adomian 1994), variational iteration method (He 1999a), Homotopy perturbation method (He 1999b), fractional differential transform method (Arigoklu and Ozkol 2007), Homotopy analysis method (Zurigat et al. 2010), new iterative method (Varsha and Hossein 2006) were proposed to obtain accurate and stable numerical solutions. Recently there has been an increasing interest in the development of numerical methods using orthogonal functions like block pulse functions, hat functions, Chebyshev polynomials, Muntz polynomials, Laguerre polynomials, Jacobi polynomials, Legendre polynomials, Bernstein polynomials, etc. (Li and Sun 2011; Tripathi et al. 2013; Doha et al. 2011, 2012; Saadatmandi and Dehghan 2010; Saadatmandi 2014; Esmaeili et al. 2011; Bhrawy and Taha 2012). Those numerical algorithms possess an interesting and significant property of transforming the fractional order differential equation under consideration into a system of algebraic equations. Once the fractional order differential equation has been converted into the system of algebraic equations, the problem of approximating the solution of fractional order differential equation becomes the task of solving the system of algebraic equations which is much easier and computationally faster. Therefore, basing numerical methods on orthogonal functions is beneficial. The new orthogonal hybrid functions (HFs), which are linear combinations of the orthogonal piecewise constant sample-and-hold functions and the orthogonal piecewise linear triangular functions, belong to the class of orthogonal polynomials of degree 1 and have been proposed by Deb et al. in (Deb et al. 2012) mainly for the analysis of delay and delay-free integer order systems and for the computation of numerical solutions of linear integer order differential equations (Deb et al. 2016). Even though very few applications of this new set of orthogonal hybrid functions has been reported in the literature, as the subsequent sections reveal, it possesses immense strength to be employed in the fractional calculus. The present work brings forward one of the important applications of the new orthogonal hybrid functions for solving multi-order fractional differential equation (FDE) of the form:

$$D^\alpha y(t) = \sum_{k=1}^{r} b_k(t) D^{\beta_k} y(t) + b_0(t)(y(t))^p + g(t), \ t \in [0,1], \ \alpha \in (n-1, n], \ n \in Z^+. \tag{1}$$

Here $\alpha > \beta_1 > \beta_2 > \cdots \beta_r$, $b_k(t)$, $b_0(t)$, $g(t)$ are given functions of independent variable $t$, $p$ is an arbitrary number, $y(t)$ is the unknown function and $D^\alpha y(t)$ is the Caputo fractional order derivative defined as (Caputo 1967)

$$D^\alpha f(t) = J^{n-\alpha} D^n f(t) = \frac{1}{\Gamma(m-\alpha)} \int_0^t (t-\tau)^{n-\alpha-1} f^{(n)}(\tau) d\tau, \tag{2}$$

where $J^\alpha f(t)$ is Riemann-Liouville fractional order integral (Oldham and Spanier 1974),

$$J^\alpha f(t) = \frac{1}{\Gamma(\alpha)} \int_0^t (t-\tau)^{\alpha-1} f(\tau) d\tau.$$

To accomplish the goal of numerically solving the multi-order fractional differential equation in (1) using HFs, we first derive a formula, by means of generalized one-shot operational matrices, expressing explicitly the fractional order integral of HFs in terms of HFs themselves and in the second step we make use of the derived generalized one-shot operational matrices to transform (1) into a system of algebraic equations whose solution is the approximate solution of (1). The remaining paper is structured as follows. Section 2 introduces new orthogonal hybrid functions. The generalized one-shot operational matrices are derived in Section 3. Based on the results of Section 3, a new numerical algorithm is proposed in Section 4. In Section 5, we test our numerical algorithm by applying it to a set of test problems consisting of linear and nonlinear multi-order fractional differential equations. Section 6 presents the significant inferences drawn from the present work.

## 2. Orthogonal hybrid functions (Deb et al. 2012, 2016)

Let us consider the time function, $f(t)$, of Lebesgue measure, which is defined on the interval, $[0,T)$. We split the interval into an $m$ number of subintervals using the constant step size, $h$.

$$[0,h), [h,2h), [2h,3h), \ldots\ldots, [jh,(j+1)h), \ldots\ldots, [(m-1)h, mh). \tag{3}$$

We now state an $m$-term hybrid function vector, $H_{(m)}(t)$, each of which represents the original function, $f(t)$, in the respective subinterval.

$$H_{(m)}(t) = [H_0(t) \quad H_1(t) \quad H_2(t) \quad \cdots \quad H_i(t) \quad \cdots \quad \cdots \quad H_{m-1}(t)]^{Tp}, \tag{4}$$

where '$Tp$' indicates transpose.

The $i^{th}$ hybrid function, $H_i(t)$, is defined as

$$H_i(t) = c_i S_i(t) + d_i T_i(t), \quad i \in [0, m-1], \tag{5}$$

where $c_i$, $d_i$ are arbitrary constants, $S_i(t)$ is the $i^{th}$ sample-and-hold function (SHF) described as

$$S_i(t) = \begin{cases} 1, & \text{if } t \in [ih, (i+1)h), \\ 0, & \text{otherwise,} \end{cases} \quad i \in [0, m-1], \tag{6}$$

and $T_i(t)$ the $i^{th}$ right-handed triangular function (TF)

$$T_i(t) = \begin{cases} \dfrac{t - ih}{h}, & \text{if } t \in [ih, (i+1)h), \\ 0, & \text{otherwise,} \end{cases} \quad i \in [0, m-1]. \tag{7}$$

Using (5) to (7), the actual function, $f(t)$, can be approximated piecewise as shown below.

$$f(t) \approx \sum_{i=0}^{m-1} H_i(t) = \sum_{i=0}^{m-1} (c_i S_i(t) + d_i T_i(t)) = C_S^T S_{(m)}(t) + C_T^T T_{(m)}(t), \tag{8}$$

where $C_S^T = [c_0 \quad c_1 \quad \cdots \quad \cdots \quad c_{m-1}]$, $C_T^T = [d_0 \quad d_1 \quad \cdots \quad \cdots \quad d_{m-1}]$, $d_i = c_{i+1} - c_i$, $c_i = f(ih)$,

$S_{(m)}(t) = [S_o(t) \quad S_1(t) \quad \cdots \quad \cdots \quad S_{m-1}(t)]^{Tp}$, $T_{(m)}(t) = [T_0(t) \quad T_1(t) \quad \cdots \quad \cdots \quad T_{m-1}(t)]^{Tp}$.

The first order integral of $S_{(m)}(t)$ is estimated in the orthogonal HF domain as follows.

$$\int_0^t S_{(m)}(\tau)d\tau = P_{1ss(m)}S_{(m)}(t) + P_{1st(m)}T_{(m)}(t), \qquad (9)$$

where $P_{1ss(m)} = h[\![0 \ 1 \ 1 \ \cdots \ \cdots \ 1]\!]_{m\times m}$, $P_{1st(m)} = h[\![1 \ 0 \ 0 \ \cdots \ \cdots \ 0]\!]_{m\times m}$,

$$[\![a \ b \ c \ d]\!] = \begin{bmatrix} a & b & c & d \\ 0 & a & b & c \\ 0 & 0 & a & b \\ 0 & 0 & 0 & a \end{bmatrix}.$$

In the same way, the right-handed triangular function vector, $T_{(m)}(t)$ is integrated once and the result is expanded into orthogonal HF series

$$\int_0^t T_{(m)}(\tau)d\tau = P_{1ts(m)}S_{(m)}(t) + P_{1tt(m)}T_{(m)}(t), \qquad (10)$$

where $P_{1ts(m)} = \dfrac{h}{2}[\![0 \ 1 \ 1 \ \cdots \ \cdots \ 1]\!]_{m\times m}$, $P_{1tt(m)} = \dfrac{h}{2}[\![1 \ 0 \ 0 \ \cdots \ \cdots \ 0]\!]_{m\times m}$.

Employing (9) and (10),

$$\begin{aligned}\int_0^t f(\tau)d\tau &\approx \int_0^t \left(C_S^T S_{(m)}(\tau) + C_T^T T_{(m)}(\tau)\right)d\tau = C_S^T \int_0^t S_{(m)}(\tau)d\tau + C_T^T \int_0^t T_{(m)}(\tau)d\tau, \\ &= C_S^T \left(P_{1ss(m)}S_{(m)}(t) + P_{1st(m)}T_{(m)}(t)\right) + C_T^T \left(P_{1ts(m)}S_{(m)}(t) + P_{1tt(m)}T_{(m)}(t)\right), \\ &= \left(C_S^T P_{1ss(m)} + C_T^T P_{1ts(m)}\right)S_{(m)}(t) + \left(C_S^T P_{1st(m)} + C_T^T P_{1tt(m)}\right)T_{(m)}(t),\end{aligned} \qquad (11)$$

where $P_{1ss(m)}$, $P_{1st(m)}$, $P_{1ts(m)}$ and $P_{1tt(m)}$ are called the complementary one-shot operational matrices of first order integration and acting as a first order integrator in the orthogonal HF domain. Orthogonal and a few operational properties of HFs are listed in Appendix 1.

3. **Generalized one-shot operational matrices for fractional integration**.

The main objective of this section is to derive an approximation using orthogonal HFs for the Riemann-Liouville fractional integral of order $\alpha$ of $f(t)$.

**Theorem 3.1** Let $\alpha \in (n-1, n]$, $n$ is an integer, and $t \in [0, T]$. The fractional integral of order $\alpha$ of the SHF vector, $S_{(m)}(t)$, is obtained in the orthogonal HF domain as

$$J^\alpha S_{(m)}(t) = \frac{1}{\Gamma(\alpha)}\int_0^t (t-\tau)^{\alpha-1} S_{(m)}(\tau)d\tau = P_{\alpha ss(m)}S_{(m)}(t) + P_{\alpha st(m)}T_{(m)}(t), \qquad (12)$$

where $P_{\alpha ss(m)} = \dfrac{h^\alpha}{\Gamma(\alpha+1)}[\![0 \ \varsigma_1 \ \varsigma_2 \ \varsigma_3 \ \cdots \ \varsigma_{m-1}]\!]$, $\varsigma_k = (k^\alpha - (k-1)^\alpha)$, $k \in [1, m-1]$,

$$P_{\alpha st(m)} = \frac{h^{\alpha}}{\Gamma(\alpha+1)}[\![1 \ \xi_1 \ \xi_2 \ \xi_3 \ \cdots \ \xi_{m-1}]\!], \ \xi_k = (k+1)^{\alpha} - 2k^{\alpha} + (k-1)^{\alpha}, \ k \in [1, m-1].$$

*Proof*

The proof is given in Appendix 2.

**Remark**. If $\alpha = 1$, then $\varsigma_k = 1$, $\forall k \in [1, m-1]$ and $\xi_k = 0$, $\forall k \in [1, m-1]$, thus, Equation (12) produces the same result as Equation (9).

**Theorem 3.2** The fractional integral of order $\alpha$ ($\alpha \in (n-1, n], n \in Z^+$) of the right-handed triangular function vector, $T_{(m)}(t)$, is expressed by means of orthogonal HFs

$$J^{\alpha}T_{(m)}(t) = \frac{1}{\Gamma(\alpha)}\int_0^t (t-\tau)^{\alpha-1} T_{(m)}(\tau)d\tau = P_{\alpha ts(m)}S_{(m)}(t) + P_{\alpha tt(m)}T_{(m)}(t), \tag{13}$$

$$P_{\alpha ts(m)} = \frac{h^{\alpha}}{\Gamma(\alpha+2)}[\![0 \ \phi_1 \ \phi_2 \ \phi_3 \ \cdots \ \phi_{m-1}]\!], \ P_{\alpha tt(m)} = \frac{h^{\alpha}}{\Gamma(\alpha+2)}[\![1 \ \psi_1 \ \psi_2 \ \psi_3 \ \cdots \ \psi_{m-1}]\!],$$

$$\phi_k = k^{\alpha+1} - (k-1)^{\alpha}(k+\alpha), \ \psi_k = (k+1)^{\alpha+1} - (k+1+\alpha)k^{\alpha} - k^{\alpha+1} + (k+\alpha)(k-1)^{\alpha}, \ k \in [1, m-1].$$

*Proof*

Appendix 3 presents the proof of **Theorem 3.2**.

**Remark**. Since $\phi_k = 1$ and $\psi_k = 0$ for $\alpha = 1$ and $k \in [1, m-1]$, the result of Equation (10) can be recovered from Equation (13).

**Theorem 3.3** The HF approximation for the Riemann-Liouville fractional order integral of the function $f(t)$ is

$$\frac{1}{\Gamma(\alpha)}\int_0^t (t-\tau)^{\alpha-1} f(\tau)d\tau \approx \left(C_S^T P_{\alpha ss(m)} + C_T^T P_{\alpha ts(m)}\right) S_{(m)}(t) + \left(C_S^T P_{\alpha st(m)} + C_T^T P_{\alpha tt(m)}\right) T_{(m)}(t). \tag{14}$$

*Proof*

Using Equation (8) in the definition of the Riemann-Liouville fractional order integral of $f(t)$ is

$$\begin{aligned} J^{\alpha}f(t) &= \frac{1}{\Gamma(\alpha)}\int_0^t (t-\tau)^{\alpha-1} f(\tau)d\tau \approx \frac{1}{\Gamma(\alpha)}\int_0^t (t-\tau)^{\alpha-1} \left(C_S^T S_{(m)}(\tau) + C_T^T T_{(m)}(\tau)\right) d\tau, \\ &= C_S^T \left(\frac{1}{\Gamma(\alpha)}\int_0^t (t-\tau)^{\alpha-1} S_{(m)}(\tau) d\tau\right) + C_T^T \left(\frac{1}{\Gamma(\alpha)}\int_0^t (t-\tau)^{\alpha-1} T_{(m)}(\tau) d\tau\right). \end{aligned} \tag{15}$$

From **Theorems 3.1** and **3.2**, we can get the expression in (14). This completes the proof.

**Remark**. The generalized HF approximation in (14) yields the HF estimate for the first order integration of $f(t)$ in special case $\alpha = 1$.

**Example 3.1**

We now test the viability of the generalized one-shot operational matrices based HF approximation ($J_{HF}^{\alpha} f(t)$) for the Riemann-Liouville fractional order integral of $f(t)$. Let us consider a time function $f(t) = t$, $t \in [0,1]$ and the step size of 0.125.

The exact fractional order integral of $f(t)$ is

$$J^{\alpha} f(t) = \frac{\Gamma(1+1)t^{1+\alpha}}{\Gamma(2+\alpha)}, \quad \alpha \in (0,5]. \tag{16}$$

Table 1 presents the $\infty$-norm of the error between the exact fractional integral and its approximation in the orthogonal HF domain for various of $\alpha$. It is noticed from Table 1 that the formula in (14) is able to approximate both the integer order and the non-integer order integral of $f(t)$ with high accuracy. Even if the HF estimate for $J^{\alpha} f(t)$ may seem like complex, it requires very less CPU usage. The approach of estimating the Riemann-Liouville fractional integral by the use of the generalized one-shot operational matrices is reliable, exact and computationally effective.

**Table 1** Performance of HF approximation

| $\alpha$ | $\left\| J^{\alpha} f(t) - J_{HF}^{\alpha} f(t) \right\|_{\infty}$ | CPU time (s) |
| --- | --- | --- |
| 0.5 | 1.11022302462516e-16 | 0.021720 |
| 1 | 0 | 0.020514 |
| 1.5 | 5.551115123125783e-17 | 0.022866 |
| 2 | 2.775557561562891e-17 | 0.021033 |
| 2.5 | 1.387778780781446e-17 | 0.023839 |
| 3 | 6.938893903907228e-18 | 0.021494 |
| 3.5 | 8.673617379884035e-19 | 0.021404 |
| 4 | 0 | 0.022049 |
| 4.5 | 1.084202172485504e-19 | 0.021662 |
| 5 | 2.168404344971009e-19 | 0.021855 |

## 4. The algorithm for solving multi-order FDEs (HFM)

We now explain the procedure of finding the approximate solution of the multi-order FDE by the application of the generalized one-shot operational matrices derived in the previous section.

Let us recall the general form of the multi-order FDE.

$$D^{\alpha} y(t) = \sum_{k=1}^{r} b_k(t) D^{\beta_k} y(t) + b_0(t)(y(t))^p + g(t), \quad t \in [0,1], \quad \alpha \in (n-1, n], \quad n \in Z^+. \tag{17}$$

For the sake of simplicity, we consider the homogeneous initial conditions $y^{(s)}(0) = 0$, $s = 0,1,2,\ldots,n-1$.

Let

$$D^{\alpha} y(t) \approx C_S^T S_{(m)}(t) + C_T^T T_{(m)}(t) \tag{18}$$

$$\begin{aligned}D^{\beta_k} y(t) &\approx J^{\alpha-\beta_k}\left(C_S^T S_{(m)}(t) + C_T^T T_{(m)}(t)\right) \\&= C_S^T\left(P_{(\alpha-\beta_k)ss(m)} S_{(m)}(t) + P_{(\alpha-\beta_k)st(m)} T_{(m)}(t)\right) + C_T^T\left(P_{(\alpha-\beta_k)ts(m)} S_{(m)}(t) + P_{(\alpha-\beta_k)tt(m)} T_{(m)}(t)\right) \\&= \left(C_S^T P_{(\alpha-\beta_k)ss(m)} + C_T^T P_{(\alpha-\beta_k)ts(m)}\right) S_{(m)}(t) + \left(C_S^T P_{(\alpha-\beta_k)st(m)} + C_T^T P_{(\alpha-\beta_k)tt(m)}\right) T_{(m)}(t)\end{aligned} \tag{19}$$

$$(y(t))^p \approx \tilde{C}_S^T S_{(m)}(t) + \tilde{C}_T^T T_{(m)}(t) \text{ (using (A6))}. \tag{20}$$

$$g(t) \approx C_{S0}^T S_{(m)}(t) + C_{T0}^T T_{(m)}(t). \tag{21}$$

$$\begin{aligned}y(t) &\approx J^{\alpha}\left(C_S^T S_{(m)}(t) + C_T^T T_{(m)}(t)\right) \\&= \left(C_S^T P_{\alpha ss(m)} + C_T^T P_{\alpha ts(m)}\right) S_{(m)}(t) + \left(C_S^T P_{\alpha st(m)} + C_T^T P_{\alpha tt(m)}\right) T_{(m)}(t)\end{aligned} \tag{22}$$

Equation (17) becomes,

$$\begin{aligned}C_S^T S_{(m)}(t) + C_T^T T_{(m)}(t) &= \sum_{k=1}^{r} b_k(t)\left(C_S^T P_{(\alpha-\beta_k)ss(m)} + C_T^T P_{(\alpha-\beta_k)ts(m)}\right) S_{(m)}(t) + C_{S0}^T S_{(m)}(t) + C_{T0}^T T_{(m)}(t) + \\&\quad \sum_{k=1}^{r} b_k(t)\left(C_S^T P_{(\alpha-\beta_k)st(m)} + C_T^T P_{(\alpha-\beta_k)tt(m)}\right) T_{(m)}(t) + b_0(t)\left(\tilde{C}_S^T S_{(m)}(t) + \tilde{C}_T^T T_{(m)}(t)\right)\end{aligned} \tag{23}$$

Estimating the variable coefficients $b_k(t)$ and $b_0(t)$ in the HF domain,

$$b_k(t) \approx C_{S1}^T S_{(m)}(t) + C_{T1}^T T_{(m)}(t),\ b_0(t) \approx C_{S2}^T S_{(m)}(t) + C_{T2}^T T_{(m)}(t). \tag{24}$$

From (23) and (24),

$$\begin{aligned}C_S^T S_{(m)}(t) + C_T^T T_{(m)}(t) &= \sum_{k=1}^{r}\left(C_{S1}^T S_{(m)}(t) + C_{T1}^T T_{(m)}(t)\right)\left(C_S^T P_{(\alpha-\beta_k)ss(m)} + C_T^T P_{(\alpha-\beta_k)ts(m)}\right) S_{(m)}(t) + \\&\quad C_{S0}^T S_{(m)}(t) + C_{T0}^T T_{(m)}(t) + \sum_{k=1}^{r}\left(C_{S1}^T S_{(m)}(t) + C_{T1}^T T_{(m)}(t)\right)\left(C_S^T P_{(\alpha-\beta_k)st(m)} + C_T^T P_{(\alpha-\beta_k)tt(m)}\right) T_{(m)}(t) + \\&\quad \left(C_{S2}^T S_{(m)}(t) + C_{T2}^T T_{(m)}(t)\right)\left(\tilde{C}_S^T S_{(m)}(t) + \tilde{C}_T^T T_{(m)}(t)\right)\end{aligned} \tag{25}$$

Employing (A5),

$$C_S^T S_{(m)}(t) + C_T^T T_{(m)}(t) = \sum_{k=1}^{r}\left(C_1 S_{(m)}(t) + (C_2 + C3) T_{(m)}(t)\right) + \left(C_{S0}^T + C_4\right) S_{(m)}(t) + \left(C_{T0}^T + C_5\right) T_{(m)}(t), \tag{26}$$

where $C_1 = C_{S1}^T .*\left(C_S^T P_{(\alpha-\beta_k)ss(m)} + C_T^T P_{(\alpha-\beta_k)ts(m)}\right)$, $C_2 = C_{T1}^T .*\left(C_S^T P_{(\alpha-\beta_k)ss(m)} + C_T^T P_{(\alpha-\beta_k)ts(m)}\right)$,

$$C_3 = \left(C_{S1}^T + C_{T1}^T\right) * \left(C_S^T P_{(\alpha-\beta_k)st(m)} + C_T^T P_{(\alpha-\beta_k)tt(m)}\right), C_5 = \left(C_{S2}^T + C_{T2}^T\right) * \tilde{C}_T^T + C_{T2}^T .* \tilde{C}_S^T,$$

$C_4 = C_{S2}^T .* \tilde{C}_S^T$ and the operator .* indicates element-wise multiplication.

Equating the coefficients of the SHF vector, $S_{(m)}(t)$, and the TF vector, $T_{(m)}(t)$,

$$C_S^T = \sum_{k=1}^{r} C_1 + \left(C_{S0}^T + C_4\right). \tag{27}$$

$$C_T^T = \sum_{k=1}^{r} (C_2 + C3) + \left(C_{T0}^T + C_5\right). \tag{28}$$

From (22), we get the piecewise linear approximation for $y(t)$ in the orthogonal HF domain.

Solving the system of algebraic equations in (27) and (28) produces the approximation for $y(t)$ in the orthogonal HF domain. As the formula in (14) works well for arbitrary order $\alpha$, the developed numerical algorithm is generic in the sense that it can be applied to classical (integer) and fractional (non-integer) multi-order differential equations.

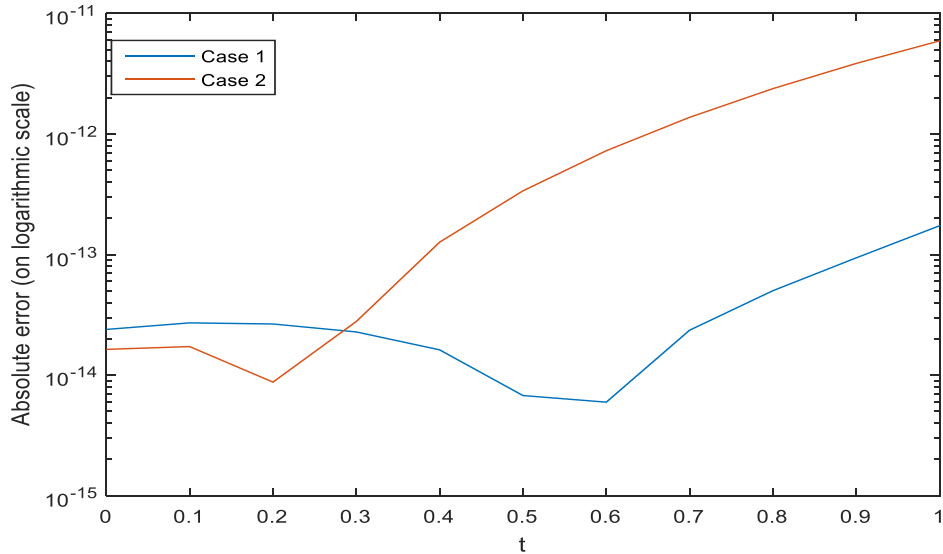

**Figure 1** Absolute error of Example 5.1.

## 5. Numerical examples

In this section, we implement the numerical algorithm developed in the former section on linear and nonlinear multi-order FDEs with constant and variable coefficients.

### Example 5.1

Consider the following linear multi-order FDE (Ford and Connolly 2009)

$$D^\alpha y(t) + D^\beta y(t) + y(t) = f(t), \ y(0) = 0, \ y'(0) = 0. \tag{29}$$

The exact solution of (27) is $y(t) = t^3$.

We consider the following three cases.

*Case 1*

$$\alpha = 2, \ \beta = 0.5, \ f(t) = t^3 + 6t + \left(\frac{3.2}{\Gamma(0.5)}\right) t^{2.5}. \tag{30}$$

The multi-order FDE in (29) with above parameters is solved using our numerical algorithm (HFM) with the step size of 0.1. The absolute error between the exact solution and the piecewise linear HF solution is computed and plotted in Figure 1. The $\infty$-norm of the absolute error is compared in Table 2 with the maximal absolute error produced in (Ford and Connolly 2009) [25] (see Example 1 and Table 1 in (Ford and Connolly 2009)), (Shiralashetti and Deshi 2016) (see Example 1 and Table 2 in (Shiralashetti and Deshi 2016)) and (Hesameddini et al. 2016) (see Example 5.2 and Table 1 in (Hesameddini et al. 2016)). It is evident in Tables 2 and 3 that our numerical algorithm is able to give more accurate approximate solution even with small number of subintervals $m = 10$ i.e. $h = 0.1$ than those achieved in (Ford and Connolly 2009), (Shiralashetti and Deshi 2016) and (Hesameddini et al. 2016).

*Case 2*

$$\alpha = 2, \ \beta = 0.75, \ f(t) = t^3 + 6t + \frac{8.533333333}{\Gamma(0.25)} t^{2.25}. \tag{31}$$

The maximal absolute error attained by our numerical algorithm (HFM1) is compared (Table 4) with that got in (Shiralashetti and Deshi 2016) (see Example 3 and Table 6 in (Shiralashetti and Deshi 2016)), (Ford and Connolly 2009) (see Example 3 and Table 3 in (Ford and Connolly 2009)) and (Bhrawy et al. 2013) (see Example 5 and Tables 6, 7 in (Bhrawy et al. 2013)). In this case too, our numerical algorithm provides superior result with the step size of 0.1. The time elapsed during the computation of HF solution of (29) (case 2) is noted and tabulated in Table 3. Figure 1 shows the absolute error (on logarithmic scale) obtained via HFM in case 2.

**Table 2** Comparison of $\infty$-norm of HFM with other methods

| Method | Step size, $h$ | Maximal absolute error |
|---|---|---|
| **HFM** | **1/10** | **1.7341e-13** |
| HWCM (Shiralashetti and Deshi 2016) | 1/512 | 1.8626e-09 |
| Method 1a (Ford and Connolly 2009) | 1/512 | 2.96e-04 |
| Method 1b (Ford and Connolly 2009) | 1/512 | 2.71e-04 |
| Method 2 (Ford and Connolly 2009) | 1/512 | 1.79e-05 |
| Method 3 (Ford and Connolly 2009) | 1/512 | 2.96e-04 |
| Method 1a(2) (Ford and Connolly 2009) | 1/512 | 8.14e-07 |
| Method 3(2) (Ford and Connolly 2009) | 1/512 | 8.14e-07 |
| RVIM* (Hesameddini et al. 2016) | - | 8.55e−10 |

*In (Hesameddini et al. 2016), 7th iterate is considered to get the approximate solution.

**Table 3** Computational time needed by HFM

| Example | Step size $h$ | CPU time (in seconds) | |
|---|---|---|---|
| | | Case 1 | Case 2 |
| 5.1 | 1/10 | 0.189338 | 0.196911 |
| 5.2 | 1/500 | 19.621723 | 20.527074 |
| 5.3 | 1/500 | 19.278202 | 15.573774 |

**Table 4** Performance of HFM in comparison with other methods for case 2 of Example 5.1

| Method | $h$ | Maximal absolute error | Method | $\alpha$ | Maximal absolute error |
|---|---|---|---|---|---|
| **HFM** | **1/10** | **5.91726e-12** | GLT(GQ)[a] (Bhrawy et al. 2013) ($N=64$) | 0 | 2.16e-07 |
| HWCM (Shiralashetti and Deshi 2016) | 1/512 | 1.8624e-09 | | 1 | 2.51e-06 |
| | | | | 2 | 7.29e-06 |
| | | | | 3 | 1.43e-05 |
| Method 1a (Ford and Connolly 2009) | 1/512 | 3.54e-03 | | | |
| Method 1b (Ford and Connolly 2009) | 1/512 | 6.93e-05 | GLT(GRQ)[b] (Bhrawy et al. 2013) ($N=64$) | 0 | 3.08e-07 |
| | | | | 1 | 4.95e-06 |
| Method 2 (Ford and Connolly 2009) | 1/512 | 1.18e-04 | | 2 | 1.80e-05 |
| | | | | 3 | 4.29e-05 |
| Method 3 (Ford and Connolly 2009) | 1/512 | 5.43e-04 | | | |
| Method 1a(2) (Ford and Connolly 2009) | 1/512 | 3.10e-06 | | | |
| Method 3(2) (Ford and Connolly 2009) | 1/512 | 5.07e-06 | | | |

[a,b] where $N$ is the degree of generalized Laguerre polynomial.

**Table 5** $\infty$-norm via HFM, ET and ER for Example 5.2

| Method | Step size $h$ | Maximal absolute error | |
|---|---|---|---|
| | | Case 1 | Case 2 |
| **HFM** | **1/500** | **2.6650e-05** | **7.006286e-06** |
| ET | 1/1000 | 0.0009980642 | 0.0009958541 |
| ER | 1/1000 | 0.0009538565 | 0.0009805249 |

**Example 5.2**

The linear multi-order FDE is

$$aD^\alpha x(t) + bD^{\alpha_2} x(t) + cD^{\alpha_1} x(t) + ex(t) = \frac{2b}{\Gamma(3-\alpha_2)} t^{2-\alpha_2} + \frac{2c}{\Gamma(3-\alpha_1)} t^{2-\alpha_1} - \frac{c}{\Gamma(2-\alpha_1)} t^{1-\alpha_1} + e(t^2 - t),$$

$x(0) = 0$, $x'(0) = -1$, $x''(0) = 2$, $x'''(0) = 0$. \hfill (32)

The exact solution for (32) is $x(t) = t^2 - t$.

The given multi-order FDE is solved using HFM for the following two cases.

*Case 1*

$$a = 1,\ b = 1,\ c = 1,\ e = 1,\ \alpha_1 = 0.77,\ \alpha_2 = 1.44,\ \alpha = 3.91. \tag{33}$$

*Case 2*

$$a = 1,\ b = 1,\ c = 0.5,\ e = 0.5,\ \alpha_1 = \sqrt{2}/20,\ \alpha_2 = \sqrt{2},\ \alpha = \sqrt{11}. \tag{34}$$

Figure 2 gives the absolute error obtained (respective computation times are provided in Table 3) via our numerical algorithm for two cases. In (El-Mesiry et al. 2005) (see Example 9 in (El-Mesiry et al. 2005)), the linear multi-order FDE in (32) is solved using ET (Euler's method with product trapezoidal quadrature formula) and ER (Euler's method with product rectangle rule) methods using the step size of 0.001. Comparing with the results of (El-Mesiry et al. 2005), our numerical algorithm provides better result (Table 5) with twice the step size used in (El-Mesiry et al. 2005).

**Table 6** Accuracy of HFM and PECE for Example 5.3

| Method | Step size $h$ | Maximal absolute error | |
|---|---|---|---|
| | | Case 1 | Case 2 |
| **HFM** | **1/500** | **1.841512e-07** | **1.965186e-07** |
| PECE | 1/1000 | 0.0004096262 | 0.0004379749 |

**Example 5.3**

The linear multi-order FDE is

$$aD^2 x(t) + bDx(t) + cD^{\alpha_2} x(t) + eD^{\alpha_1} x(t) + kx(t) = f(t),\ x(0) = 1,\ x'(0) = 0, \tag{35}$$

where $f(t) = a + bt + \dfrac{c}{\Gamma(3 - \alpha_2)} t^{2-\alpha_2} + \dfrac{e}{\Gamma(3 - \alpha_1)} t^{2-\alpha_1} + k(1 + 0.5t^2)$.

The given linear non-homogenous multi-order FDE has the analytical solution $x(t) = 1 + 0.5t^2$.

We consider the following two cases.

*Case 1*

$$a = 1,\ b = 3,\ c = 2,\ e = 1,\ k = 5,\ \alpha_1 = 0.0159,\ \alpha_2 = 0.1379. \tag{36}$$

*Case 2*

$$a = 0.2,\ b = 1,\ c = 1,\ e = 0.5,\ k = 2,\ \alpha_1 = 0.00196,\ \alpha_2 = 0.07621. \tag{37}$$

The accuracy of the HF approximate solution obtained in both cases (shown in Figure 3, Table 6 and the respective elapsed times are given in Table 3) is better than the accuracy of the numerical solution acquired by Predict Evaluate-Correct Evaluate method (PECE) in (EL-Sayed et al. 2004) (see Example 8 and Tables 7, 8 in (EL-Sayed et al. 2004)).

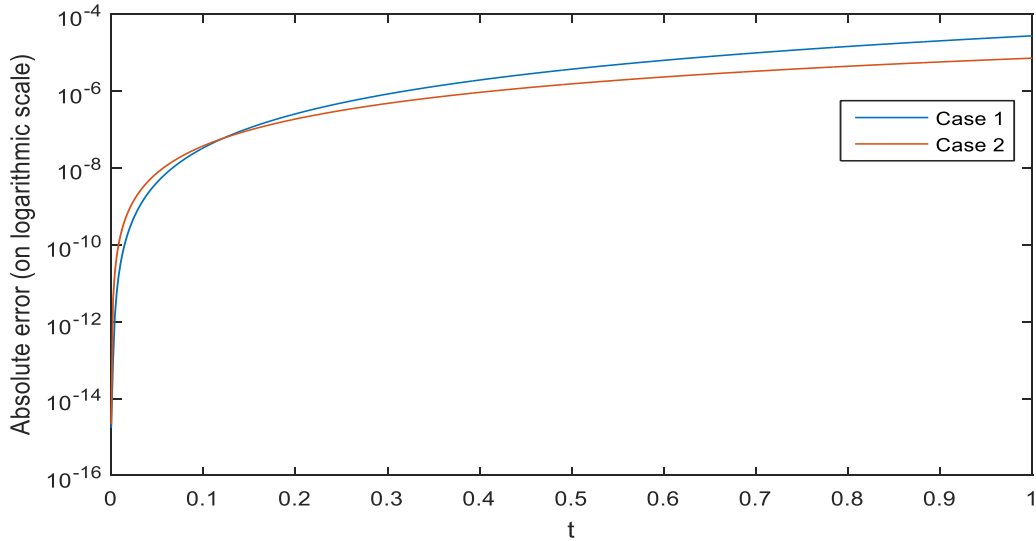

**Figure 2** Absolute error of Example 5.2.

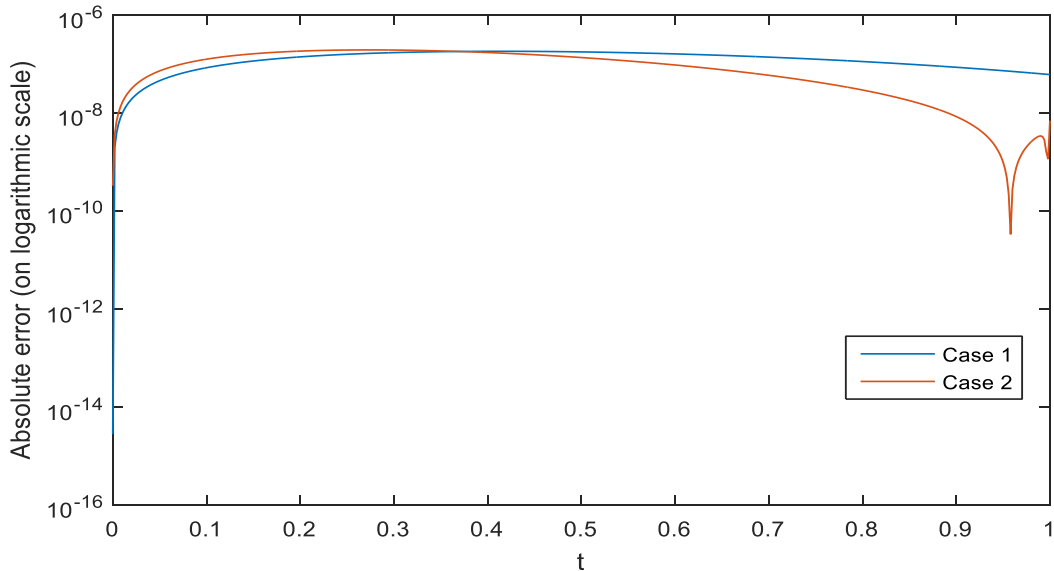

**Figure 3** Absolute error of Example 5.3.

### Example 5.4

The nonlinear multi-order FDE is

$$aD^\alpha x(t) + bD^{\alpha_2} x(t) + cD^{\alpha_1} x(t) + e(x(t))^3 = f(t),\ x(0)=0,\ x'(0)=0,\ x''(0)=0, \tag{38}$$

where $f_2(t) = \left(2at^{3-\alpha}/\Gamma(4-\alpha)\right) + \left(2b/\Gamma(4-\alpha_2)\right)t^{3-\alpha_2} + \left(2c/\Gamma(4-\alpha_1)\right)t^{3-\alpha_1} + e\left(t^3/3\right)^3$.

The exact solution for the following five cases is $x(t) = t^3/3$.

*Case 1*

$a = 1$, $b = 2$, $c = 0.5$, $e = 1$, $\alpha_1 = 0.00196$, $\alpha_2 = 0.07621$, $\alpha = 2$. (39)

*Case 2*

$a = 1$, $b = 0.1$, $c = 0.2$, $e = 0.3$, $\alpha_1 = \sqrt{5}/5$, $\alpha_2 = \sqrt{2}/2$, $\alpha = 2$. (40)

Using the step size of 1/10, the nonlinear non-homogenous multi-order FDE in (38) is solved for the above two cases and the corresponding absolute errors are compared in Figure 4(a). Table 7 presents the maximal absolute errors gained via our algorithm and those obtained by ET and ER methods in (EL-Sayed et al. 2005) (see Example 2 in (EL-Sayed et al. 2005)), through Adomian decomposition method (ADM) and modified Podlubny's numerical method (PNM) in (El-Sayed et al. 2010) (see Example 2 and Tables 3 and 6 in (El-Sayed et al. 2010)) and via a numerical method transforming the fractional order differential equations into a system of first order ordinary differential equations in (Javidi and Nyamoradi 2013) (see Example 1 and Tables 1 and 3 in (Javidi and Nyamoradi 2013)). In terms of accuracy and computation time (Table 8), the performance of our numerical algorithm (HFM) is far more superior to that of ET and ER (in (EL-Sayed et al. 2005)), PNM and ADM (in (El-Sayed et al. 2010)) and NM (in (Javidi and Nyamoradi 2013)).

**Table 7** Error analysis of Example 5.4 in case 1 and case 2.

| Method | Step size $h$ | Maximal absolute error | |
|---|---|---|---|
| | | Case 1 | Case 2 |
| **HFM** | **1/10** | **7.205347e-14** | **5.268008e-14** |
| ET | 1/1000 | 0.0008924007 | 0.0009717941 |
| ER | 1/1000 | 0.0007891357 | 0.0009438396 |
| PNM | 1/2000 | 0.000399235 | 0.000388881 |
| ADM* | - | 0.000150218 | 5.74351e-06 |
| NM | 1/2000 | 9.39e−5 | 2.6866e−4 |

*The series solution obtained by ADM is truncated to $N = 3$.

**Table 8** CPU time taken by HFM and other methods.

| Method | CPU time (in seconds) | |
|---|---|---|
| | Case 1 | Case 2 |
| **HFM** | **0.195641** | **0.199136** |
| PNM | 894.999 | 952.907 |
| ADM | 478.577 | 506.954 |
| NM | 5.9380 | 5.9060 |

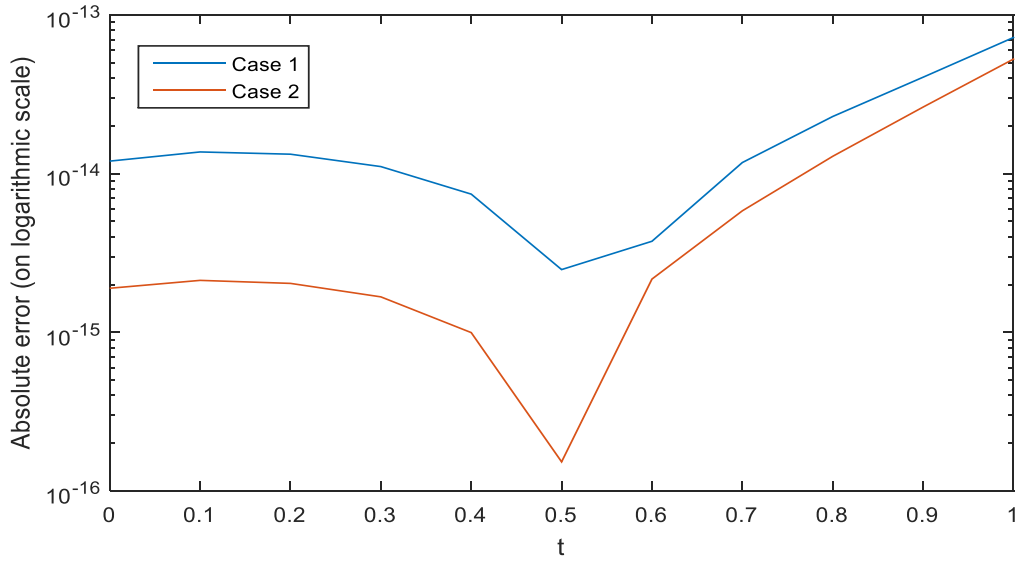

**Figure 4(a)** Absolute error of Example 5.4 in case 1 and case 2.

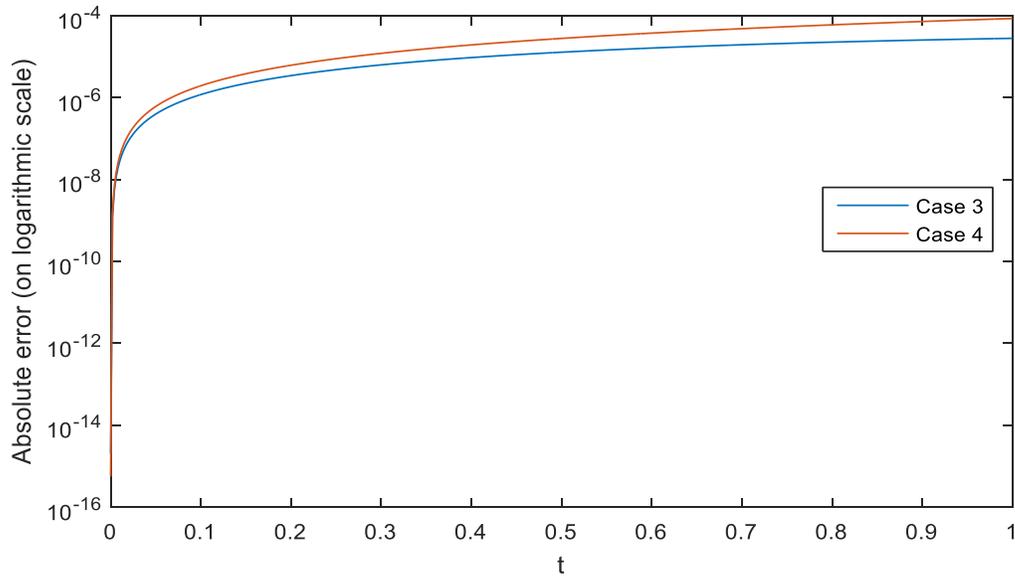

**Figure 4(b)** Absolute error of Example 5.4 in case 3 and case 4.

*Case 3*

$a = 1$, $b = 2$, $c = 0.5$, $e = 1$, $\alpha_1 = 0.00196$, $\alpha_2 = 1.07621$, $\alpha = 2.55..$ (41)

*Case 4*

$a = 1$, $b = 0.1$, $c = 0.2$, $e = 0.3$, $\alpha_1 = \sqrt{7}/7$, $\alpha_2 = \sqrt{7}/2$, $\alpha = \sqrt{7}$. (42)

Table 9 and Figure 4(b) emphasis that our numerical algorithm has better convergence than 2E and 3E methods in (EL-Sayed et al. 2004) (see Example 10 in (EL-Sayed et al. 2005)). The computation time needed by HFM in case 3 and case 4 are 19.655548 and 16.739572 seconds, respectively.

**Table 9** Error analysis of Example 5.4 in case 3 and case 4.

| Method | Step size h | Maximal absolute error | |
|---|---|---|---|
| | | Case 3 | Case 4 |
| **HFM** | 1/500 | **2.80956210e-05** | **8.5593407e-05** |
| 2E | 1/1000 | 0.00134739300 | 0.00150859400 |
| 3E | 1/1000 | 0.00120052700 | 0.00149777500 |

**Table 10** Error analysis of Example 5.5.

| Method | Step size h | Maximal absolute error | |
|---|---|---|---|
| | | Case 1 | Case 2 |
| **HFM** | 1/300 | **4.96926379e-06** | **1.62141126e-07** |
| PECE | 1/1000 | 0.000006325524 | 0.00057616830 |

**Example 5.5**

Consider the following nonlinear multi-order FDE with variable coefficients

$$aD^2 x(t) + bD^{\alpha_2} x(t) + c\left(D^{\alpha_1} x(t)\right)^2 + e(x(t))^3 = 2at + \frac{2bt^{3-\alpha_2}}{\Gamma(4-\alpha_2)} + c\left(\frac{2t^{3-\alpha_1}}{\Gamma(4-\alpha_1)}\right)^2 + e\left(\frac{t^3}{3}\right)^3, \quad (43)$$

subject to the initial conditions $x(0) = 0$, $x'(0) = 0$.

The exact solution is $x(t) = t^3/3$.

The piecewise linear HF approximate solution of (43) is obtained by using the numerical algorithm for the following cases.

*Case 1*

$a = 1$, $b = 1$, $c = 1$, $e = 1$, $\alpha_1 = 0.555$, $\alpha_2 = 1.455$. \hfill (44)

*Case 2*

$a = 1$, $b = 0.5$, $c = 0.5$, $e = 0.5$, $\alpha_1 = 0.276$, $\alpha_2 = 1.999$. \hfill (45)

As shown in Table 10 and Figure 5, our method gave more accurate results than those results by Predictor Evaluate-Corrector Evaluate method (PECE) in (EL-Sayed et al. 2004) (see Example 6 and Tables 3 and 4 in (EL-Sayed et al. 2004)). The time elapsed during computation of HF solution in case 1 and case 2 are given in Table 11.

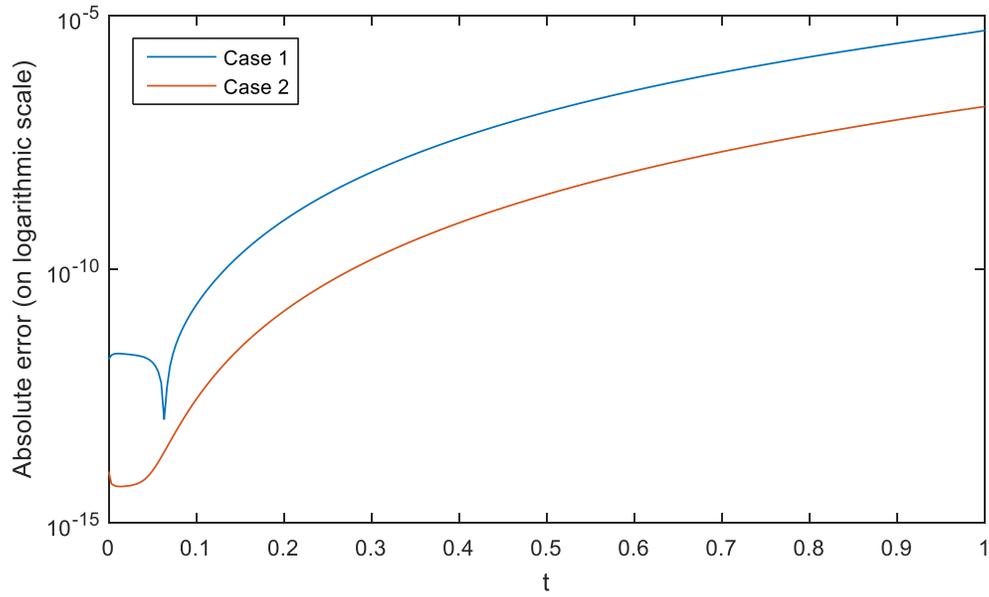

**Figure 5** Absolute error of Example 5.5.

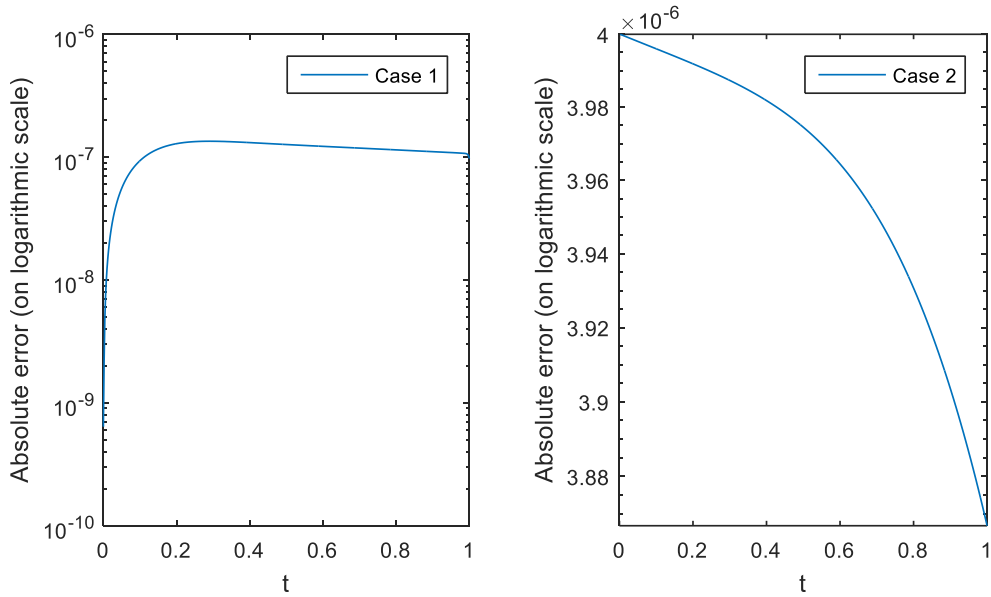

**Figure 6** Error analysis of Example 5.6.

**Table 11** CPU time taken by HFM for Example 5.5 and 5.6

| Example | Step size $h$ | CPU time (in seconds) | |
|---|---|---|---|
| | | Case 1 | Case 2 |
| 5.5 | **1/300** | 6.813894 | 5.074030 |
| 5.6 | 1/500 | 36.837632 | 48.637273 |

**Table 12** Error analysis of Example 5.6

| Method | Step size $h$ | Maximal absolute error | |
|---|---|---|---|
| | | Case 1 | Case 2 |
| **HFM** | **1/500** | **1.34738857e-07** | **4.00000005e-06** |
| ET (Example 3 in EL-Sayed et al. 2005) | 1/1000 | 0.00031542780 | 0.00050354 |
| ER (Example 3 in EL-Sayed et al. 2005) | 1/1000 | 0.00002789852 | 0.0004844666 |

**Example 5.6**

The linear multi-order FDE with variable coefficients is

$$aD^2 x(t) + b(t)Dx(t) + c(t)D^{\alpha_2} x(t) + e(t)D^{\alpha_1} x(t) + k(t)x(t) = f(t), \quad x(0) = 2, \quad x'(0) = 0, \tag{46}$$

where $f(t) = -a - b(t)t - \dfrac{c(t)t^{2-\alpha_2}}{\Gamma(3-\alpha_2)} - \dfrac{e(t)t^{2-\alpha_1}}{\Gamma(3-\alpha_1)} + k(t)(2 - 0.5t^2)$.

The given problem possesses the closed form solution $x(t) = 2 - 0.5t^2$.

Figure 6 displays the absolute error of HF solution for the following cases.

*Case 1*

$$a = 0.1, \ b(t) = t, \ c(t) = 1 + t, \ e(t) = t^2, \ k(t) = (1+t)^2 \ \alpha_1 = 0.781, \ \alpha_2 = 0.891. \tag{47}$$

*Case 2*

$$a = 5, \ b(t) = \sqrt{t}, \ c(t) = t^2 - t, \ e(t) = 3t, \ k(t) = t^3 - t, \ \alpha_1 = \sqrt{7}/70, \ \alpha_2 = \sqrt{13}/13. \tag{48}$$

The results in Figure 6 and Tables 11 and 12 ensure that the developed numerical algorithm is suitable to handle a wide variety of multi-order fractional differential equations.

### 6. Conclusions

We have the following concluding remarks.

- Tables 1 to 11 prove that, from viewpoint of solving numerically multi-order FDEs, the proposed numerical algorithm, despite the complexity of the problem under consideration, can provide more accurate results than Adomian decomposition method and Podlubny's numerical method (El-Sayed et al. 2010), reconstruction of variational iteration method (Hesameddini et al. 2016), Predictor Evaluate-Correct Evaluate method (PECE) (EL-Sayed et al. 2004), Euler's method with product trapezoidal quadrature formula) and ER (Euler's method with product rectangle rule (EL-Sayed et al. 2005), generalized Laguerre-spectral methods (Bhrawy et al. 2013), systems-based decomposition schemes (Ford and Connolly 2009) and Haar wavelet collocation method (Shiralashetti and Deshi 2016).

- Though the present work has been carried out from the standpoint of pure mathematics, the important results (the generalized one-shot hybrid function operational matrices and the numerical algorithm) of the work may be worth using for precise analysis of physical processes modelled by multi-order FDEs.
- The successful utilization of the new orthogonal hybrid functions for solving a special class of fractional order differential equations (i.e. multi-order FDEs) encourages to explore the unexplored applications areas (fractional order system analysis, fractional order system identification, fractional order controller design, etc.) of the new orthogonal hybrid functions.

**Appendix 1** Basic properties of HFs.

The components of SHF vector, $S_{(m)}(t)$, and TF vector, $T_{(m)}(t)$, have orthogonal property.

$$\int_0^T S_i(t)S_j(t)dt = \begin{cases} h, & \text{if } i == j, \\ 0, & \text{otherwise}, \end{cases} \int_0^T T_i(t)T_j(t)dt = \begin{cases} \dfrac{h}{3}, & \text{if } i == j, \\ \dfrac{h}{6}, & \text{otherwise}, \end{cases} i \in [0, m-1]. \tag{A1}$$

The product, $S_i(t)S_j(t)$, where $i, j \in [0, m-1]$, is expressed via HFs

$$S_i(t)S_j(t) = \begin{cases} S_i(t), & \text{if } i == j, \\ 0, & \text{otherwise}. \end{cases} \tag{A2}$$

Similarly,

$$T_i(t)T_j(t) = \begin{cases} T_i(t), & \text{if } i == j, \\ 0, & \text{otherwise}, \end{cases} i, j \in [0, m-1]. \tag{A3}$$

Whereas the product, $S_i(t)T_j(t)$, is estimated in the orthogonal HF domain as

$$S_i(t)T_j(t) = \begin{cases} T_j(t), & \text{if } i == j, \\ 0, & \text{othwersie}, \end{cases} i, j \in [0, m-1]. \tag{A4}$$

The product of two functions $h(t) = f_1(t)f_2(t)$ can be approximated by HFs as given in the next equation.

$$\begin{aligned} f_1(t)f_2(t) &\approx \left(C_{S1}^T S_{(m)}(t) + C_{T1}^T T_{(m)}(t)\right)\left(C_{S2}^T S_{(m)}(t) + C_{T2}^T T_{(m)}(t)\right) \\ &= \left(C_{S1}^T \cdot * C_{S2}^T\right) S_{(m)}(t) + \left(C_{S1}^T \cdot * C_{T2}^T + C_{T1}^T \cdot * C_{S2}^T + C_{T1}^T \cdot * C_{T2}^T\right) T_{(m)}(t) \end{aligned} \tag{A5}$$

The $n^{th}$ power of function, $g(t)$, ($g(t) \in C[0,T]$) is expanded into orthogonal HFs using the following expression.

$$(g(t))^n \approx C_S^T S_{(m)}(t) + C_T^T T_{(m)}(t), \tag{A6}$$

where $C_S^T = [c_0 \quad c_1 \quad \cdots \quad \cdots \quad c_{m-1}]$, $C_T^T = [d_0 \quad d_1 \quad \cdots \quad \cdots \quad d_{m-1}]$, $c_i = (g(ih))^n$.

**Appendix 2** Proof of Theorem 3.1

The fractional integral of order $\alpha$ of $S_0(t)$ is

$$J^\alpha S_0(t) = \frac{1}{\Gamma(\alpha)} \int_0^t (t-\tau)^{\alpha-1} S_0(\tau) d\tau = \begin{cases} 0, & \text{for } t = 0 \text{ and } j = 0 \\ \frac{h^\alpha}{\Gamma(\alpha+1)} (j^\alpha - (j-1)^\alpha), & \text{for } t > 0, j > 0. \end{cases} \tag{A7}$$

Evaluating the expression in (A7) at $j = 1, 2, 3, \ldots, m-1$ yields the following coefficients.

$$c_0 = 0, \quad c_1 = \frac{h^\alpha}{\Gamma(\alpha+1)}, \quad c_2 = \frac{h^\alpha}{\Gamma(\alpha+1)} (2^\alpha - 1^\alpha), \tag{A8}$$

$$c_j = \frac{h^\alpha}{\Gamma(\alpha+1)} (j^\alpha - (j-1)^\alpha), \quad j = 3, 4, \ldots, m-1. \tag{A9}$$

The difference between the consecutive coefficients,

$$d_0 = c_1 - c_0 = \frac{h^\alpha}{\Gamma(\alpha+1)}, \quad d_1 = c_2 - c_1 = \frac{h^\alpha}{\Gamma(\alpha+1)} (2^\alpha - 1^\alpha) - \frac{h^\alpha}{\Gamma(\alpha+1)}, \tag{A10}$$

$$d_2 = c_3 - c_2 = \frac{h^\alpha}{\Gamma(\alpha+1)} (3^\alpha - 2^\alpha) - \frac{h^\alpha}{\Gamma(\alpha+1)} (2^\alpha - 1^\alpha), \tag{A11}$$

$$d_j = c_{j+1} - c_j = \frac{h^\alpha}{\Gamma(\alpha+1)} ((j+1)^\alpha - 2j^\alpha + (j-1)^\alpha), \quad j = 3, 4, \ldots, m-1. \tag{A12}$$

We can approximate $J^\alpha S_0(t)$ in terms of HFs,

$$J^\alpha S_0(t) = [c_0 \quad c_1 \quad c_2 \quad \cdots \quad c_{m-1}] S_{(m)}(t) + [d_0 \quad d_1 \quad d_2 \quad \cdots \quad d_{m-1}] T_{(m)}(t). \tag{A13}$$

Substituting the expressions for $c_i$ and $d_i$ in (A13),

$$J^\alpha S_0(t) = \frac{h^\alpha}{\Gamma(\alpha+1)} [0 \quad 1 \quad (2^\alpha - 1) \quad \cdots \quad (j^\alpha - (j-1)^\alpha) \quad \cdots \quad ((m-1)^\alpha - (m-2)^\alpha)] S_{(m)}(t) +$$

$$\frac{h^\alpha}{\Gamma(\alpha+1)} [1 \quad (2^\alpha - 2) \quad \cdots \quad ((j+1)^\alpha - 2j^\alpha + (j-1)^\alpha) \quad \cdots \quad ((m)^\alpha - 2(m-1)^\alpha + (m-2)^\alpha)] T_{(m)}(t)$$

$$\tag{A14}$$

Rewriting (A14),

$$J^\alpha S_0(t) = \frac{h^\alpha}{\Gamma(\alpha+1)}[0 \quad \varsigma_1 \quad \varsigma_2 \quad \varsigma_3 \quad \cdots \quad \varsigma_{m-1}]S_{(m)}(t) + \frac{h^\alpha}{\Gamma(\alpha+1)}[1 \quad \xi_1 \quad \xi_2 \quad \xi_3 \quad \cdots \quad \xi_{m-1}]T_{(m)}(t)$$

where $\varsigma_k = (k^\alpha - (k-1)^\alpha)$, $\xi_k = (k+1)^\alpha - 2k^\alpha + (k-1)^\alpha$, $k \in [1, m-1]$. (A15)

Carrying out fractional integration on the remaining terms and expressing the results via orthogonal HFs,

$$J^\alpha S_1(t) = \frac{h^\alpha}{\Gamma(\alpha+1)}[0 \quad 0 \quad \varsigma_1 \quad \varsigma_2 \quad \cdots \quad \varsigma_{m-2}]S_{(m)}(t) + \frac{h^\alpha}{\Gamma(\alpha+1)}[0 \quad 1 \quad \xi_1 \quad \xi_2 \quad \cdots \quad \xi_{m-2}]T_{(m)}(t).$$

(A16)

$$\vdots$$

$$J^\alpha S_{m-2}(t) = \frac{h^\alpha}{\Gamma(\alpha+1)}[0 \quad 0 \quad \cdots \quad \cdots \quad 0 \quad \varsigma_1]S_{(m)}(t) + \frac{h^\alpha}{\Gamma(\alpha+1)}[0 \quad \cdots \quad \cdots \quad 0 \quad 1 \quad \xi_1]T_{(m)}(t). \quad (A17)$$

$$J^\alpha S_{m-1}(t) = \frac{h^\alpha}{\Gamma(\alpha+1)}[0 \quad 0 \quad \cdots \quad \cdots \quad 0 \quad 0]S_{(m)}(t) + \frac{h^\alpha}{\Gamma(\alpha+1)}[0 \quad \cdots \quad \cdots \quad 0 \quad 0 \quad 1]T_{(m)}(t). \quad (A18)$$

Therefore,

$$J^\alpha S_{(m)}(t) = [J^\alpha S_0(t) \quad J^\alpha S_2(t) \quad J^\alpha S_3(t) \quad \cdots \quad J^\alpha S_{m-1}(t)]^{Tp} = P_{\alpha ss(m)}S_{(m)}(t) + P_{\alpha st(m)}T_{(m)}(t), (A19)$$

where

$$P_{\alpha ss(m)} = \frac{h^\alpha}{\Gamma(\alpha+1)}[[0 \quad \varsigma_1 \quad \varsigma_2 \quad \varsigma_3 \quad \cdots \quad \varsigma_{m-1}]], \quad P_{\alpha st(m)} = \frac{h^\alpha}{\Gamma(\alpha+1)}[[1 \quad \xi_1 \quad \xi_2 \quad \xi_3 \quad \cdots \quad \xi_{m-1}]].$$

This proves **Theorem 3.1**.

**Appendix 3** Proof of Theorem 3.2

We get the following expression upon performing fractional integration on $T_0(t)$,

$$J^\alpha T_0(t) = \frac{1}{\Gamma(\alpha)}\int_0^t (t-\tau)^{\alpha-1} T_0(\tau) d\tau = \begin{cases} 0, & \text{for } t = 0, \\ \frac{h^\alpha}{\Gamma(\alpha+2)}(j^{\alpha+1} - (j-1)^\alpha(j+\alpha)), & \text{for } t > 0, j \in [1, m-1]. \end{cases}$$

(A20)

In the orthogonal HF domain, $J^\alpha T_0(t)$ is expressed as

$$J^\alpha T_0(t) = [c_0 \quad c_1 \quad c_2 \quad \cdots \quad c_{m-1}]S_{(m)}(t) + [d_0 \quad d_1 \quad d_2 \quad \cdots \quad d_{m-1}]T_{(m)}(t), \quad (A21)$$

where $c_j = \dfrac{h^\alpha}{\Gamma(\alpha+2)}\left(j^{\alpha+1} - (j-1)^\alpha(j+\alpha)\right)$, $j = 3, 4, \ldots, m-1$,

$$d_j = c_{j+1} - c_j = \dfrac{h^\alpha}{\Gamma(\alpha+2)}\left((j+1)^{\alpha+1} - (j+1+\alpha)j^\alpha - j^{\alpha+1} + (j+\alpha)(j-1)^\alpha\right).$$

Using the expressions for $c_j$ and $d_j$ and rearranging,

$$J^\alpha T_0(t) = \Psi\left[0 \quad 1 \quad (2^{\alpha+1} - (2+\alpha)) \quad \cdots \quad (j^{\alpha+1} - (j-1)^\alpha(j+\alpha)) \quad \cdots \quad \Omega_1\right] S_{(m)}(t) +$$
$$\Psi\left[1 \quad (2^{\alpha+1} - (3+\alpha)) \quad \cdots \quad (j+1)^{\alpha+1} - (j+1+\alpha)j^\alpha - j^{\alpha+1} + (j+\alpha)(j-1)^\alpha \quad \cdots \quad \Omega_2\right] T_{(m)}(t)$$

where $\Psi = \dfrac{h^\alpha}{\Gamma(\alpha+2)}$, $\Omega_1 = \left((m-1)^{\alpha+1} - (m-2)^\alpha(m-1+\alpha)\right)$,

$$\Omega_2 = m^{\alpha+1} - (m+\alpha)(m-1)^\alpha - (m-1)^{\alpha+1} + (m-1+\alpha)(m-2)^\alpha. \tag{A22}$$

We now rewrite (A22),

$$J^\alpha T_0(t) = \dfrac{h^\alpha}{\Gamma(\alpha+2)}\left[0 \quad \phi_1 \quad \phi_2 \quad \phi_3 \quad \cdots \quad \phi_{m-1}\right] S_{(m)}(t) + \dfrac{h^\alpha}{\Gamma(\alpha+2)}\left[1 \quad \psi_1 \quad \psi_2 \quad \psi_3 \quad \cdots \quad \psi_{m-1}\right] T_{(m)}(t)$$

where $\phi_k = k^{\alpha+1} - (k-1)^\alpha(k+\alpha)$, $\psi_k = (k+1)^{\alpha+1} - (k+1+\alpha)k^\alpha - k^{\alpha+1} + (k+\alpha)(k-1)^\alpha$. (A23)

Following the same procedure, the remaining components of $T_{(m)}(t)$ can be fractional integrated and the resulting expressions can be approximated via HFs.

The fractional integral of order $\alpha$ of $T_{(m)}(t)$ in HFs domain is

$$J^\alpha T_{(m)}(t) = P_{\alpha ts(m)} S_{(m)}(t) + P_{\alpha tt(m)} T_{(m)}(t), \tag{A24}$$

where

$$P_{\alpha ts(m)} = \dfrac{h^\alpha}{\Gamma(\alpha+2)}\left[\!\left[0 \quad \phi_1 \quad \phi_2 \quad \phi_3 \quad \cdots \quad \phi_{m-1}\right]\!\right], \quad P_{\alpha tt(m)} = \dfrac{h^\alpha}{\Gamma(\alpha+2)}\left[\!\left[1 \quad \psi_1 \quad \psi_2 \quad \psi_3 \quad \cdots \quad \psi_{m-1}\right]\!\right].$$

This complete the proof.

**References**

[1] Adomian, G.: Solving Frontier Problems of Physics: The Decomposition Method. Kluwer Academic Publishers, Boston (1994)

[2] Arigoklu, A., Ozkol, I.: Solution of fractional differential equations by using differential transform method. Chaos Soliton. Fract. 34, 1473–1481 (2007)

[3] Bhrawy, A.H., Taha, T.M.: An operational matrix of fractional integration of the Laguerre polynomials and its application on a semi-infinite interval. Mathematical Sciences 6:41 (2012)

[4] Bhrawy, A.H., Baleanu, D., Assas, L.M.: Efficient generalized Laguerre-spectral methods for solving multi-term fractional differential equations on the half line. Journal of Vibration and Control 0(0), 1–13 (2013)

[5] Caputo, M.: Linear models of dissipation whose Q is almost frequency. Part II. J Roy Austral Soc. 13, 529–539 (1967)

[6] Deb, A., Ganguly, A., Sarkar, G., Biswas, A.: Numerical solution of third order linear differential equations using generalized one-shot operational matrices in orthogonal hybrid function domain. Appl. Math. Comput. 219, 1485-1514 (2012)

[7] Deb, A., Roychoudhury, S., Sarkar, G.: Analysis and Identification of Time-Invariant Systems, Time-Varying Systems, and Multi-Delay Systems Using Orthogonal Hybrid Functions. Theory and Algorithms with MATLAB. Springer International Publishing AG, Switzerland (2016)

[8] Diethelm, K., Ford, N.J, Freed, A.D.: A predictor-corrector approach for the numerical solution of fractional differential equations. Nonlinear Dyn. 29, 3–22 (2002)

[9] Doha, E.H., Bhrawy, A.H., Ezz-Eldien, S.S.: Efficient Chebyshev spectral methods for solving multi-term fractional orders differential equations. Applied Mathematical Modelling 35, 5662-5672 (2011)

[10] Doha, E.H., Bhrawy, A.H., Ezz-Eldien, S.S.: A new Jacobi operational matrix: An application for solving fractional differential equations. Applied Mathematical Modelling 36, 4931-4943 (2012)

[11] EL-Sayed, M.A., EL-Mesiry, A.E.M., EL-Saka, H.A.A.: Numerical solution for multi-term fractional (arbitrary) orders differential equations. Computational and Applied Mathematics 23(1), 33–54 (2004)

[12] El-Mesiry, A.E.M., El-Sayed, A.M.A., El-Saka, H.A.A.: Numerical methods for multi-term fractional (arbitrary) orders differential equations. Applied Mathematics and Computation 160, 683–699 (2005)

[13] El-Sayed, A.M.A., El-Kalla, I.L., Ziada, E.A.A.: Analytical and numerical solutions of multi-term nonlinear fractional orders differential equations. Applied Numerical Mathematics 60, 788–797 (2010)

[14] Esmaeili, S., Shamsi, M., Luchko, Y.: Numerical solution of fractional differential equations with a collocation method based on Muntz polynomials. Computers and Mathematics with Applications 62, 918-929 (2011)

[15] Ford, N.J., Connolly, J.A.: Systems-based decomposition schemes for the approximate solution of multi-term fractional differential equations. Journal of Computational and Applied Mathematics 229, 382–391 (2009)

[16] Gaul, L., Klein, P., Kempfle, S.: Impulse response function of an oscillator with fractional derivative in damping description. Mech. Res. Comm. 16 (5), 4447–4472 (1989)

[17] Gaul, L., Klein, P., Kempfle, S.: Damping description involving fractional operators. Mech. Systems Signal Process. 5 (2), 8–88 (1991)

[18] He, J.H.: Variational iteration method – a kind of non-linear analytical technique: Some examples. Int. J. Non-Linear Mech. 34, 699-708 (1999a)

[19] He, J.H.: Homotopy perturbation technique. Comput. Method Appl. Mech. Eng. 178, 257-262 (1999b)

[20] Hesameddini, E., Rahimi, A., Asadollahifard, E.: On the convergence of a new reliable algorithm for solving multi-order fractional differential equations. Commun. Nonlinear Sci. Numer. Simulat. 34, 154–164 (2016)

[21] Javidi, M., Nyamoradi, N.: A numerical scheme for solving multi-term fractional differential equations. Commun. Frac. Calc. 4 (1), 38–49 (2013)

[22] Li, Y., Sun, N.: Numerical solution of fractional differential equations using the generalized block pulse operational matrix. Computers and Mathematics with Applications 62, 1046–1054 (2011)

[23] Odibat, Z.M., Momani, S.: An algorithm for the numerical solution of differential equations of fractional order. J. Appl. Math. Informatics 26, 15–27 (2008)

[24] Oldham, K.B., Spanier, J.: The Fractional Calculus: Theory and Applications of Differentiation and Integration to Arbitrary Order. Dover Publications, New York (1974)

[25] Podlubny, I.: Fractional Differential Equations. Academic Press, New York (1999)

[26] Saadatmandi, A., Dehghan, M.: A new operational matrix for solving fractional-order differential equations. Computers and Mathematics with Applications 59, 1326-1336 (2010)

[27] Saadatmandi, A.: Bernstein operational matrix of fractional derivatives and its applications. Applied Mathematical Modelling 38, 1365-1372 (2014)

[28] Shiralashetti, S.C., Deshi, A.B.: An efficient Haar wavelet collocation method for the numerical solution of multi-term fractional differential equations. Nonlinear Dyn. 83, 293–303 (2016)

[29] Suarez, L.E., Shokooh, A.: An eigenvector expansion method for the solution of motion containing fractional derivatives. ASME. J. Appl. Mech. 64, 629–635 (1997)

[30] Tripathi, M.P., Baranwal, V.K., Pandey, R.K., Singh, O.P.: A new numerical algorithm to solve fractional differential equations based on operational matrix of generalized hat functions. Commun Nonlinear. Sci. Numer. Simulat. 18, 1327–1340 (2013)

[31] Varsha, D-G., Hossein, J.: An iterative method for solving nonlinear functional equations. J. Math. Anal. Appl. 316, 753–763 (2006)

[32] Zurigat, M., Momani, S., Odibat, Z., Alawneh, A.: The homotopy analysis method for handling systems of fractional differential equations. Appl. Math. Modell. 34, 24–35 (2010)